\documentclass[11pt,leqno]{article}
\begin{document}


\def\nexto{\kern -0.54em}

\newcommand{\R}{{\sf I\hspace{-.15em}R}}
\newcommand{\N}{{\sf I\hspace{-.15em}N}}
\newcommand{\C}{{\sf I\hspace{-.5em}C}}
\newcommand{\BB}{{\sf I\hspace{-.15em}B}}
\def\SS{{\sf I\hspace{-.4em}S}}

\newcommand{\Z}{\bf Z}

\def\T{{\bf T}}
\def\diam{\mbox{\rm diam}}
\def\rr{{\cal R}}
\def\mt{{\Lambda}}
\def\e{\emptyset}
\def\sd{{\bf S}^{d-1}}
\def\so{{\bf S}^1}
\def\dQ{\partial Q}
\def\dk{\partial K}
\def\endofproof{{\rule{6pt}{6pt}}}
\def\ts{\tilde{\sigma}}
\def\tr{\tilde{r}}
\def\sigt{\ts^{r}}
\def\di{\displaystyle}
\def\dist{\mbox{\rm dist}}
\def\tU{\widetilde{U}}
\def\hU{\widehat{U}}
\def\tS{\tilde{S}}
\def\tP{\widetilde{\Pi}}
\def\sa+{\Sigma_A^+}
\def\u-{\overline{u}}
\def\du{\frac{\partial}{\partial u}}
\def\dv{\frac{\partial}{\partial v}}
\def\dt{\frac{d}{d t}}
\def\dx{\frac{\partial}{\partial x}}
\def\z-{\overline{z}}
\def\te{\tilde{e}}
\def\tv{\tilde{v}}
\def\tu{\tilde{u}}
\def\tw{\tilde{w}}
\def\con{\mbox{\rm const }}
\def\Box{\spadesuit}
\def\nn{{\cal N}}
\def\mm{{\cal M}}
\def\kk{{\cal K}}
\def\ll{{\cal L}}
\def\vv{{\cal V}}
\def\nab{\nn^{(J,\mu)}_{ab}}
\def\mab{\mm_{ab}}
\def\ma{\mm_{ab}}
\def\lab{L_{ab}}
\def\tlab{\tilde{L}_{ab}}
\def\mabn{\mm_{ab}^N}
\def\man{\mm_a^N}
\def\labn{L_{ab}^N}
\def\hZ{\widehat{Z}}
\def\tz{\tilde{Z}}
\def\hz{\hat{z}}
\def\hg{\hat{\gamma}}
\def\mmu{m_{\mu,b,J}}
\def\tF{\tilde{F}}
\def\tf{\tilde{f}}
\def\tp{\tilde{p}}
\def\ff{{\cal F}}
\def\i{{\bf i}}
\def\jj{{\bf j}}
\def\ttau{\tilde{\tau}}
\def\tt{{\cal T}}
\def\uu{{\cal U}}
\def\hu{\widehat{U}}
\def\wloc{W_{\epsilon}}
\def\pp{{\cal P}}
\def\tg{\tilde{\gamma}}  
\def\aa{{\cal A}}
\def\tV{\widetilde{V}}
\def\cc{{\cal C}}
\def\tC{\widetilde{\cc}}
\def\Deo{\dd_\epsilon(\Omega)}
\def\sdk{S_{\dk}(\Omega)}
\def\lae{\Lambda_{\epsilon}}
\def\ep{\epsilon}
\def\tr{\tilde{R}}
\def\oo{{\cal O}}

\def\be{\begin{equation}}
\def\ee{\end{equation}}
\def\beqn{\begin{eqnarray*}}
\def\eeqn{\end{eqnarray*}}

\def\cc{{\cal C}}
\def\gi{\gamma^{(i)}}
\def\ii{{\imath }}
\def\jj{{\jmath }}
\def\tc{\tilde{C}}
\def\II{{\cal I}}
\def\ccij{ \cc_{i'_0,j'_0}[\eta]}
\def\hc{\hat{\cc}}
\def\dd{{\mathcal D}}
\def\la{\langle}
\def\ra{\rangle}
\def\bs{\bigskip}
\def\xio{\xi^{(0)}}
\def\xo{x^{(0)}}
\def\zo{z^{(0)}}
\def\do{\partial \Omega}
\def\dk{\partial K}
\def\dl{\partial L}
\def\ll{{\cal L}}
\def\kk{{\cal K}}
\def\rr{{\cal R}}
\def\cK{\hat{K}}
\def\kk{{\cal K}}
\def\pr{{\rm pr}}
\def\ff{{\cal F}}
\def\G{{\cal G}}
\def\dist{{\rm dist}}
\def\dds{\frac{d}{ds}}
\def\con{{\rm const}\;}
\def\Con{{\rm Const}\;}
\def\di{\displaystyle}
\def\oo{{\cal O}}
\def\hess{\mbox{\rm Hess}}
\def\endofproof{{\rule{6pt}{6pt}}}
\def\vm{\varphi^{(m)}}
\def\km{k^{(m)}}
\def\dm{d^{(m)}}
\def\kam{\kappa^{(m)}}
\def\dem{\delta^{(m)}}
\def\xim{\xi^{(m)}}
\def\ep{\epsilon}
\def\tt{\tilde{t}}
\def\tx{\tilde{x}}
\def\tp{\tilde{p}}
\def\eep{\hat{\epsilon}}
\def\ep{\epsilon}
\def\hu{\hat{u}}
\def\tS{\widetilde{S}}
\def\tt{{\cal T}}
\def\ts{\tilde{\sigma}}

\def\ms{\medskip}
\def\sn{{\SS}^{n-1}}
\def\sN{{\SS}^{N-1}}
\def\be{\begin{equation}}
\def\ee{\end{equation}}
\def\beqn{\begin{eqnarray}}
\def\eeqn{\end{eqnarray}}
\def\beqn*{\begin{eqnarray*}}
\def\eeqn*{\end{eqnarray*}}
\def\endofproof{{\rule{6pt}{6pt}}}
\def\nv{\nabla \varphi}
\def\wuloc{W^u_{\mbox{\footnotesize\rm loc}}}
\def\wsloc{W^s_{\mbox{\footnotesize\rm loc}}}

\def\dist{\mbox{\rm dist}}
\def\diam{\mbox{\rm diam}}
\def\pr{\mbox{\rm pr}}
\def\supp{\mbox{\rm supp}}
\def\Arg{\mbox{\rm Arg}}
\def\In{\mbox{\rm Int}}
\def\Im{\mbox{\rm Im}}
\def\span{\mbox{\rm span}}
\def\spec{\mbox{\rm spec}\,}
\def\Re{\mbox{\rm Re}}
\def\var{\mbox{\rm var}}
\def\conf{\mbox{\footnotesize\rm const}}
\def\Conf{\mbox{\footnotesize\rm Const}}
\def\Lip{\mbox{\rm Lip}}
\def\con{\mbox{\rm const}\;}
\def\li{\mbox{\rm li}} 
\def\ex{\mbox{\rm extd}}

\def\saa{\Sigma_A^+}
\def\sa{\Sigma_A}
\def\san{\Sigma^-_A}

\def\vxij{\varphi_{\xi,j}}
\def\tz{\tilde{z}}
\def\txi{\tilde{\xi}}
\def\tphi{\tilde{\phi}}
\def\tPhi{\widetilde{\Phi}}
\def\tmt{\tilde{\Lambda}}
\def\tOm{\widetilde{\Omega}}
\def\tal{\tilde{\alpha}}

\def\Xr{X^{(r)}}
\def\hs{\hat{\sigma}}
\def\hl{\hat{l}}
\def\tnu{\tilde{\nu}}

\def\kmax{\kappa_{\max}}
\def\kmin{\kappa_{\min}}
\def\ecc{\mbox{\rm ecc}}
\def\tB{\widetilde{B}}
\def\hh{{\cal H}}
\def\thh{\widetilde{\cal H}}
\def\hE{\widehat{E}}

\def\naf{\nabla f(z)}
\def\so{\sigma_0}
 \def\tp{\tilde{p}}
\def\hcij{\hat{\cc}_{i,j}}
\def\Xo{X^{(0)}}
\def\z1{z^{(1)}}
\def\Vo{V^{(0)}}
\def\Yo{Y{(0)}}
\def\tPsi{\tilde{\Psi}}
\def\hX{\hat{X}}
\def\hx{\hat{x}}

\def\pr{\mbox{\rm pr}}
\def\dK{\partial K}
\def\tK{\tilde{K}}
\def\tpsi{\tilde{\psi}}
\def\tq{\tilde{q}}
\def\tsi{\tilde{\sigma}}
\def\iii{\sf i}
\def\clip{C^{\mbox{\footnotesize \rm Lip}}}
\def\Lip{\mbox{\rm Lip}}
\def\lip{\mbox{\footnotesize\rm Lip}}
\def\span{\mbox{ \rm span}}
\def\beqn{\begin{eqnarray*}}
\def\eeqn{\end{eqnarray*}}
\def\yo{y^{(0)}}

\begin{center}
{\Large\bf  Non-integrability of open billiard flows\\ and Dolgopyat type estimates} 
\end{center}

\begin{center}
{\sc Luchezar Stoyanov}\\
{\it University of Western Australia, Crawley WA 6009, Australia\\
(e-mail: stoyanov@maths.uwa.edu.au)}
\end{center}


\bigskip

\noindent
{\it Abstract.} We consider open billiard flows in $\R^n$  and show that  the standard symplectic form 
$d\alpha$ in $\R^n$ satisfies a specific non-integrability condition over their non-wandering sets $\mt$.  
This allows to use the main result  in \cite{kn:St3}  and obtain Dolgopyat type estimates for spectra of Ruelle transfer operators under 
simpler conditions.
 We also describe a class of open billiard flows in $\R^n$ ($n \geq 3$) satisfying a certain pinching
condition, which in turn implies that the (un)stable laminations over the non-wandering set are $C^1$.

\ms

\section{Introduction}

\renewcommand{\theequation}{\arabic{section}.\arabic{equation}}

It is well-known that hyperbolic billiard flows in compact domains (e.g. Sinai billiards 
in Euclidean spaces or on tori) are non-integrable, just like
contact Anosov flows (see e.g. \cite{kn:KB} or Appendix B in \cite{kn:L}). However, when considering contact flows over  basic sets $\mt$ 
the general non-degeneracy of the contact form (which implies the non-integrability) 
does not say much about the dynamics of the flow over $\mt$. It is much more natural, and it turns out to be important as well, to look at the 
restriction of the contact form over tangent vectors to $\mt$ (see Sect. 2 for the definition). 
This is what we do here for non-wandering sets of open billiard flows. The motivation to study this kind of non-integrability 
comes from \cite{kn:St3} which deals with spectral estimates of Ruelle transfer operators for flows on basic sets (see 
Sect. 6 below for some details).

Let $K$ be a subset of ${\R}^{n}$ ($n\geq 2$) of the form $K = K_1 \cup K_2 \cup \ldots \cup K_{k_0}$, 
where $K_i$ are compact strictly convex disjoint domains in $\R^{n}$ with $C^2$ {\it boundaries} 
$\Gamma_i = \dk_i$  and $k_0 \geq 3$. Set 
$\Omega = \overline{{\R}^n \setminus K}$ and $\Gamma = \partial K$. We assume that $K$ satisfies the 
following (no-eclipse) condition: 
$${\rm (H)} \quad \quad\qquad  
\cases{ \mbox{\rm for every pair $K_i$, $K_j$ of different connected components 
of $K$ the convex hull of}\cr \mbox{\rm $K_i\cup K_j$ has no common points with any other connected component of $K$. }\cr}$$
With this condition, the {\it billiard flow} $\phi_t$ defined on the {\it sphere bundle} $S(\Omega)$ in the standard way is called an 
{\it open billiard flow}. It has singularities, however its restriction to the {\it non-wandering set} $\Lambda$ has only 
simple discontinuities at reflection points.  Moreover, $\Lambda$  is compact, $\phi_t$ is hyperbolic and transitive on $\Lambda$, and  
it follows from \cite{kn:St1} that $\phi_t$ is  non-lattice and therefore by  a result of Bowen  \cite{kn:B}, it is topologically weak-mixing 
on $\Lambda$.

Our main aim in this paper is to show that the open billiard flow  always satisfies a 
certain non-integrability condition on $\mt$. Let $d\alpha$ be the {\it standard symplectic form} on $T(\R^n) = \R^n \times \R^n$.

\bs

\noindent
{\sc Theorem 1.1.} {\it  There exist $z_0\in \mt$ and  $\mu > 0$  such that  for any unit tangent vector $b\in E^u(z_0)$ to $\mt$ 
there exists a unit tangent vector $a\in E^s(z_0)$ to $\mt$  with $|d\alpha (a,b)| \geq  \mu$}.

\bs

If the map $\mt \ni x \mapsto E^u(x)$ is $C^1$, then  the invariance of $d\alpha$ along the flow 
implies that the points $z_0\in \mt$ with the above property form an open and dense subset of $\mt$.
Theorem 1.1 is established by means of a certain pairing of points on the strong stable and 
unstable manifolds of an appropriately chosen point $z_0$ - see Sect. 3 and Lemma 3.1 there
for details. As a consequence of this and the main result in \cite{kn:St3} one gets Dolgopyat
type spectral estimates for pinched open billiard flows -- see Sect. 6. 

It is well-known that in general the maps $\mt \ni x \mapsto E^u(x)$ (or   $E^s(x)$) are only
H\"older continuous (see e.g \cite{kn:HPS} or \cite{kn:PSW}). The following {\it pinching condition} implies stronger
regularity properties of these maps. 

\ms

\noindent
{\sc (P)}:  {\it There exist  constants $C > 0$ and $0 < \alpha \leq \beta$ such that for every $x\in \mt$ we have
$$\frac{1}{C} \, e^{\alpha_x \,t}\, \|u\| \leq \| d\phi_{t}(x)\cdot u\| \leq C\, e^{\beta_x\,t}\, \|u\|
\quad, \quad  u\in E^u(x) \:\:, t > 0 \;,$$
for some constants $\alpha_x, \beta_x > 0$ depending on $x$ but independent of $u$ and $t$ with
$\alpha \leq \alpha_x \leq \beta_x \leq \beta$ and $2\alpha_x - \beta_x \geq \alpha$ for all $x\in \mt$.}

\ms

For example in the case of contact flows $\phi_t$, it follows from the results 
in \cite{kn:Ha2} (see also \cite{kn:Ha1}) that assuming (P), the map $\mt \ni x \mapsto E^u(x)$ is 
$C^{1+\ep}$ with 
$\ep = 2 \alpha/\beta -1 > 0$ (in the  sense that this map has a linearization at any $x\in \mt$ 
that depends H\"older
continuously on $x$). The same applies to the map $\mt \ni x \mapsto E^s(x)$.

Notice that when $n = 2$ (then the local unstable manifolds are one-dimensional) this condition 
is always satisfied. 
It turns out that for $n \geq 3$ the condition (P) is always satisfied when the minimal distance between  distinct connected 
components of $K$ is relatively large compared to the maximal sectional curvature of $\dk$ 
(see Proposition 1.2 below). 
An analogue of the latter
for manifolds $M$ of strictly negative curvature would be to require that the sectional curvature is between 
$-K_0$ and $-a\, K_0$ for some constants $K_0 > 0$ and $a\in (0,1)$. It follows from the arguments in  \cite{kn:HP} that when 
$a = 1/4$  the geodesic flow on $M$ satisfies the pinching condition (P).

Set $\di d_{i,j} = \dist(K_i, K_j)$ and $\di d_{0} = \min_{i\neq j} d_{i,j}$.
Since every $K_i$ is strictly convex, the operator  
$L_x : T_x(\dk) \longrightarrow T_x(\dk)$, $L_xu = (\nabla_u \nu)(x)$,
of the {\it second fundamental form} is positive definite with respect to the outward unit normal field $\nu(y)$, 
$y\in \dk$. Then  $k(x,u) = \la L_x u, u\ra$ is the {\it normal curvature} of $\dk$ at $x$ in the direction of $u \in T_x(\dk)$, $\| u \| = 1$. 
Set
$$\kmin = \min_{x\in \dk} \:\:\min_{u\in T_x(\dk), \|u\|=1} \langle L_x(u),  u\rangle \quad , \quad
\kmax = \max_{x\in \dk} \:\:\max_{u\in T_x(\dk), \|u\|=1} \langle L_x(u),  u\rangle\;,$$
where $\langle \cdot, \cdot\rangle$ is the {\it standard inner product} in $\R^n$.

Before continuing, notice that the condition (H) implies the existence of a global constant
$\varphi_0 \in (0,\pi/2)$ such that  for any $x\in \mt$ and any reflection point $q$  of the billiard trajectory 
$\gamma(x)$ generated by $x$ the angle $\varphi$ between the reflected direction of $\gamma(x)$ at $q$ and the
outward normal to $\dk$ at $q$ satisfies $\varphi \leq \varphi_0$.
Set $\mu_0 = 2 \cos \varphi_0\, \kmin$   and   $\di \lambda_0 = \frac{1}{d_0} + \frac{2 \, \kmax }{\cos\varphi_0} \;.$

Let $a > 0$ be such that $d_{i,j} \leq d_0 + a$ for all $i,j= 1, \ldots,k_0$, $i \neq j$. Below we 
assume that  $d_0$ is large compared to $a$ and $\kmax$, so that 
\be
[1+(d_0+a)\, \lambda_0]^{d_0+a} < (1+d_0\, \mu_0)^{2d_0}\;.
\ee
(Notice that when $a = r\, d_0$, $0 < r < 1$, then the above holds for all sufficiently large $d_0$,
assuming $\kmax$ and $\kmin$ are uniformly bounded above and below, respectively, by positive constants.)

In Section 5 below we prove  the following

\bs

\noindent
{\sc Proposition 1.2.} {\it  Assume that  {\rm (1.1)} holds and the boundary $\dk$ is $C^3$.   Then the open billiard flow $\phi_t$ in the exterior of $K$
satisfies the condition {\rm (P)} on its non-wandering set $\mt$. Moreover, for any $x\in \mt$ we can choose
$\alpha_x = \alpha_0$ and $\beta_x = \beta_0$, where
$\alpha_0 =  \frac{ \ln (1+ d_0\, \mu_0)}{d_0 +a}$ and}  
$\beta_0 =   \frac{ \ln (1+ (d_0+a)\, \lambda_0)}{d_0}\;.$

\bs

This is relatively easy to derive from a formula for the growth of the differential of the flow on 
unstable manifolds (see Proposition 5.1). The latter can be proved using an argument similar to that 
in the Appendix in \cite{kn:St2} (dealing with the two-dimensional case), and also can be easily derived 
from more general facts about the evolution of unstable vectors for multidimensional dispersing
billiards (see e.g. \cite{kn:BCST}).

Section 2 below contains some basic definitions and an example which concerns the geometry of 
the non-wandering set $\mt$.  Sections 3 and 4 are devoted to the proof of Theorem 1.1.
In Section 5 we use some well known formulae of Sinai for curvature operators related to unstable 
manifolds of dispersing billiards to prove Proposition 1.2.
Section 6 deals with Dolgopayt type estimates for pinched open billiard flows -- these are straightforward
consequences of \cite{kn:St3} and  the considerations in Sect. 3 below.

\bs

\noindent
{\it Acknowledgement.} The author is grateful to the referee whose comments and criticism led
to a significant improvement of the first version of the paper.

\section{Preliminaries}
\setcounter{equation}{0}

\subsection{Basic definitions}

Let $M$ be a $C^1$ complete  Riemann manifold,  and  $\phi_t : M \longrightarrow M$ ($t\in \R$) a 
$C^1$ flow on $M$. A $\phi_t$-invariant closed subset $\mt$ of $M$ is called {\it hyperbolic} if 
$\mt$ contains
no fixed points  and there exist  constants $C > 0$ and $0 < \lambda < 1$ such that  there exists a 
$d\phi_t$-invariant decomposition 
$T_xM = E^0(x) \oplus E^u(x) \oplus E^s(x)$ of $T_xM$ ($x \in \mt$) into a direct sum of non-zero 
linear subspaces,
where $E^0(x)$ is the one-dimensional subspace determined by the direction of the flow at $x$, 
$\| d\phi_t(u)\| \leq C\, \lambda^t\, \|u\|$ 
for all  $u\in E^s(x)$ and $t\geq 0$, and $\| d\phi_t(u)\| \leq C\, \lambda^{-t}\, \|u\|$ for all 
$u\in E^u(x)$ and  $t\leq 0$.

A non-empty compact $\phi_t$-invariant hyperbolic subset $\mt$ of $M$ which is not a single 
closed orbit is called a {\it basic set} for 
$\phi_t$ if $\phi_t$ is transitive on $\mt$ and $\mt$ is locally maximal, i.e. there 
exists an open neighbourhood $V$ of
$\mt$ in $M$ such that $\mt = \cap_{t\in \R} \phi_t(V)$. 

For $x\in \Lambda$ and a sufficiently small $\epsilon > 0$ let 
$$\wloc^s(x) = \{ y\in M : d (\phi_t(x),\phi_t(y)) \leq \epsilon \: \mbox{\rm for all }
\: t \geq 0 \; , \: d (\phi_t(x),\phi_t(y)) \to_{t\to \infty} 0\: \}\; ,$$
$$\wloc^u(x) = \{ y\in M : d (\phi_t(x),\phi_t(y)) \leq \epsilon \: \mbox{\rm for all }
\: t \leq 0 \; , \: d (\phi_t(x),\phi_t(y)) \to_{t\to -\infty} 0\: \}$$
be the (strong) {\it stable} and {\it unstable manifolds} of size $\epsilon$. Then
$E^u(x) = T_x \wloc^u(x)$ and $E^s(x) = T_x \wloc^s(x)$. 
Given $z \in \mt$, let $\exp^u_z : E^u(z) \longrightarrow W^u_{\ep_0}(z)$  and
$\exp^s_z : E^s(z) \longrightarrow W^s_{\ep_0}(z)$ be the corresponding
{\it exponential maps}. A  vector $b\in E^u(z)\setminus \{ 0\}$ is called  {\it tangent to $\mt$} at
$z$ if there exist infinite sequences $\{ v^{(m)}\} \subset  E^u(z)$ and 
$\{ t_m\}\subset \R\setminus \{0\}$
such that $\exp^u_z(t_m\, v^{(m)}) \in \mt \cap W^u_{\ep}(z)$ for all $m$, $v^{(m)} \to b$ and 
$t_m \to 0$ as $m \to \infty$.  It is easy to see that a vector $b\in E^u(z)\setminus \{ 0\}$ is  
tangent to $\mt$ at
$z$ iff there exists a $C^1$ curve $z(t)$, $0\leq t \leq a$, in $W^u_{\ep}(z)$ for some $a > 0$ 
with $z(0) = z,\: \dot{z}(0) = b$, and $z(t_n) \in \mt$ for some sequence 
$\{t_n\}_{n=1}^\infty \subset (0,a]$ with $t_n \to 0$ as $n \to \infty$. 
Tangent vectors to $\mt$ in $E^s(z)$ are defined similarly.  
Denote by $\hE^u(z)$ (resp. $\hE^s(z)$) the {\it set of all 
vectors $b\in E^u(z)\setminus \{ 0\}$ (resp. $b\in E^s(z)\setminus \{ 0\}$) tangent to} $\mt$ at $z$.

\bs

\noindent
{\it Remark 1.} Although we have not sought to construct particular examples, it appears that
in general the set of unit tangent vectors to $\mt$ does not have to be closed in the bundle 
$E^u_\mt$ (or $E^s_\mt$). That is, there may exist a point  $z \in \mt$, a
sequence $\{ z_m\} \subset W^u_{\ep}(z)\cap \mt$ and for each $m$ a unit vector
$\xi_m$ tangent to $\mt$ at $z_m$ such that $z_m \to z$ and $\xi_m \to \xi$ as $m \to \infty$,
however $\xi$ is not tangent to $\mt$ at $z$. 

\bs

Next, assume that $K$ and $\Omega$ are as in Sect. 1.
The {\it non-wandering set} $\Lambda$ for the flow $\phi_t$ is the set of those $x\in S(\Omega)$ such that the trajectory 
$\{ \phi_t(x) : t\in \R\}$ is bounded. Notice that the natural projection of $\phi_t$ on the quotient space $S(\Omega)/\sim$, where
$\sim$ is the equivalence relation $(q,v) \sim (p,w)$ iff $q=p$ and $v = w$ or 
$q = p\in \dk$ and $v$ and $w$ are symmetric with respect to $T_q(\dk)$, is continuous. Moreover
whenever both $x$ and $\phi_t(x)$ are in the interior of $S(\Omega)$ and sufficiently
close to $\Lambda$, the map $y \mapsto \phi_t(y)$ is smooth on a neighbourhood of $x$.
It follows from results of Sinai (\cite{kn:Si1}, \cite{kn:Si2}) that $\Lambda$ is a hyperbolic
set for $\phi_t$, and it is easily seen that 
$\Lambda$ is the maximal compact $\phi_t$-invariant subset of $S(\Omega)$.
Moreover, it follows from the natural symbolic coding for the natural section of
the flow (the so called billiard ball map) that the periodic points are dense in $\Lambda$, 
and $\phi_t$ is transitive on $\lambda$.
Thus, $\Lambda$ is a basic set for $\phi_t$ and the classical theory of hyperbolic flows
applies in the case under consideration (see e.g. Part 4 in \cite{kn:KH}).

\subsection{An example}

Here we briefly describe a (non-trivial)  example from \cite{kn:St4} which shows that in 
general for every $z\in \mt$ the space $\span(\hE^u_{\tmt}(z))$
generated by the vectors in $E^u(z)$ tangent to $\mt$ could be a proper subspace
of $E^u(z)$.


\bs

\noindent
{\bf Example 2.1.} (\cite{kn:St4}) Assume that $n = 3$ and there exists a plane $\alpha$ 
such that each of the domains $K_j$ is symmetric with respect to $\alpha$.  Setting $K' = K\cap \alpha$ and
$\Omega' = \Omega \cap \alpha$, it is easy to observe that every billiard trajectory generated by
a point in $\mt$ is entirely contained in $\alpha$. That is, $\mt = \mt'$, where $\mt'$ is the 
non-wandering set for the open billiard flow in $\Omega'$. Thus, 
$\dim(\span(\hE^u_\mt(z))) = 1 < \dim(E^u(z)) = 2$ 
for any $z \in \mt$. This example is of course trivial, since $\mt$ is contained in the flow-invariant
submanifold $S^*(\Omega')$ of $S^*(\Omega)$.

However with a small local perturbation of the boundary $\dk$ of $K$ we can get a non-trivial example.
Choosing standard cartesian coordinates $x,y,z$ in $\R^3$, we may assume that
$\alpha$ is given by the equation $z = 0$, i.e. $\alpha = \R^2\times \{0\}$.
Let $\pr_1: S^*(\R^3) \sim \R^3 \times \SS^2 \longrightarrow \R^3$ be the natural projection, and let
$C = \pr_1(\mt)$. We may choose the coordinates $x,y$ in the plane $\alpha = \{ z=0\}$ so that
the line $y =0$ is tangent to $K'_1$ and $K'_2$ and $K'$ is contained in
the half-plane $y \geq 0$. Let $q_1\in K'_1$ and $q_2\in K'_2$ be such that $[q_1,q_2]$ is the
shortest segment connecting $K'_1$ and $K'_2$. Take a point $q_1' \in \dk'_1$ close to
$q_1$ and such that the $y$-coordinate of $y'_1$ is smaller than that of $q_1$. Consider the {\it open
arc} $\aa$ on $\dk'_1$ connecting $q_1$ and $q_1'$. It is clear that $\aa \cap C = \e$.

Let $f: \R^3\longrightarrow \R^3$ be a $C^1$ (we can make it even $C^\infty$) diffeomorphism with 
$f(x) = x$ for all $x$ outside a small open set $U$ such that $q_1 \in \overline{U}$ 
and $U \cap \dk' \subset \aa$. Then for any $q\in C$ the tangent 
planes $T_q(\dk)$ and  $T_q(\partial \tK)$ coincide. We can choose $f$ so that $\tK_i = f(K_i) = K_i$ 
for $i > 1$,  $\tK_1 = f(K_1)$ is strictly convex, and  $\tnu(f(q)) \notin \alpha$ for $q \in \aa$ arbitrarily close to $q_1$. Here 
$\tnu$ is the outward  unit normal field to $\partial \tK$. 

One can then show  that the non-wandering set $\tmt$ for the billiard
flow $\tphi_t$ in the closure $\tOm$ of the exterior of $\tK$ in $\R^3$ coincides with $\mt$ 
(\cite{kn:St4}). Thus, 
$\dim(\span(\hE^u_{\tmt}(z))) = 1 < \dim (E^u(z))$ for any $z \in \tmt$.
However, it is clear from the construction that $S^*(\alpha \cap \Omega)$ is not invariant 
with respect to the billiard flow $\tphi_t$. Moreover, it is not difficult to see that there is 
no two-dimensional submanifold $\tal$ of $\tOm$ such that $S^*(\tal)$ is $d\tphi_t$-invariant and
$\mt \subset S^*(\tal)$; see Section 4 in \cite{kn:St4} for details.

\section{Non-integrability of open billiard flows}
\setcounter{equation}{0}

In this section we prove Theorem 1.1.

Let $K \subset \R^n$ be as in Sect. 1. For any $x\in  \Gamma = \partial K$ we will denote by 
$\nu(x)$ the {\it outward unit normal} to $\Gamma$ at $x$.
Given $\delta > 0$ 
denote by $S_\delta(\Omega)$ the set of those $(x,u)\in S(\Omega)$ such
that there exist $y\in \Gamma$ and $t \geq 0$ with 
$y+t u = x$, $y+ s u \in \R^n \setminus K$ for all
$s\in (0,t)$ and $\langle u, \nu_\Gamma (y) \rangle \geq \delta$. 

\bs

\noindent
{\it Remark 2.} Notice that the condition (H) implies the existence of a constant $\delta_0 > 0$ 
depending only on the obstacle $K$ such that any $(x,u) \in S(\Omega)$ whose backward and forward 
billiard trajectories both have a common point with $\Gamma$ belongs to $S_{\delta_0}(\Omega)$.

\bs

For $\epsilon \in (0,d_0/2)$ set
$$\Deo = \{ x = (q,v) \in S(\Omega) : \dist (q,\dk) > \epsilon\} \quad , 
\quad  \mt_\epsilon = \mt\cap \Deo\;.$$
In what follows in order to avoid ambiguity and unnecessary complications we will consider stable
and unstable manifolds only for points $x$ in $\Deo$ or $\lae$; this will be enough for our purposes.

Fix for a moment arbitrary $\ep$, $\delta$  and $\lambda$ so that
 \be
0 < \delta \leq \ep < \lambda < \frac{d_0}{2}\;.
 \ee 
 We will see later how small these numbers need to be.
 
Consider an arbitrary point $\sigma_0 = (\xo , \xio) \in \mt_\ep$ such that 
$\zo= \xo + \lambda\, \xio\in \dk$,
$\xio = - \nu(\zo)$ and $\xo + t\, \xio \in \R^n \setminus K$ for all $t \in [0,\lambda)$. 
I.e. the billiard trajectory 
generated by $\sigma_0$ is perpendicular to $\dk$ at $\zo$ and so the reflected direction at 
$\zo$ is $-\xio$. Notice that  there exist
such points\footnote{In fact, it is not difficult to see that the union of the orbits of 
such points $\sigma_0$
is a dense subset of $\mt$.}, e.g. we can take $\xo$ on the shortest segment between two connected 
components $K_i$ and $K_j$ ($i \neq j$) with $\xio$ parallel to that segment.
The local submanifolds $U = W^u_\delta (\sigma_0)$ and $S = \wloc^s(\sigma_0)$ have the form
$$U = \{ (x,\nu_X(x)) : x\in X\} \:\:\:, \:\:\:  S = \{ (y,\nu_Y(y) ) : y\in Y\}$$
for some smooth local $(n-1)$-dimensional  submanifolds $X$ and $Y$ in $\R^n$, where $\nu_X$ and $\nu_Y$
are continuous unit normal fields on $X$ and $Y$. Moreover the {\it second fundamental form} $L^{(X)}$ of
$X$ with respect to $\nu_X$ (resp. $L^{(Y)}$ of $Y$ with respect to $\nu_Y$) is positive 
(resp. negative) definite. Finally, we have $\xo \in X\cap Y$, $X$ and $Y$ are tangent at $\xo$ and 
$\nu_X(\xo) = \nu_Y(\xo) = \xio$. Since the tangent planes to $X$ and $Y$ at $x_0$ are parallel to the
tangent plane to $\dk$ at $\zo$, we have $T_{\xo}X = T_{\xo}Y = T_{\zo}(\dk)$.

Consider the {\it inversion} $\iii : S(\Omega) \longrightarrow S(\Omega)$ defined 
by ${\iii}(x,\xi) = (x,-\xi)$.
It follows from the general properties of stable (unstable) manifolds that for any  $(x,\xi)\in S$
(or $U$) sufficiently close to $\sigma_0$ we have $\iii\circ \phi_{2\lambda}(x,\xi)\in U$ ($S$, 
respectively).  In other words the shift along the billiard flow  $\phi_t$ of the convex front $S$
along the normal field $\nu_Y(y)$ in $2\lambda$ units coincides locally with the inversion of the
convex front $X$ near  $\xo$. Using Sinai's formula (\cite{kn:Si2}, cf. also \cite{kn:SiCh}) in the 
particular situation considered here we have
$$L_{\xo}^{(X)}(u) = \frac{B_{\xo}(u)}{I + \lambda\, B_{\xo}(u)}\;,$$
where
$$B_{\xo}(u) = \frac{L_{\xo}^{(Y)}(u)}{1+ \lambda \,L_{\xo}^{(Y)}(u)} + 2 L_{\zo}(u) \:\:\:\:, \:\:\: 
u\in T_{\xo}X\;.$$
It is well-known (see \cite{kn:Si2}) that the curvature operators of strong unstable manifolds of 
$\phi_t$ are uniformly bounded, so there exists a global constant $C > 0$ such that
$2 L_{\zo}(u) \leq B_{\xo}(u) \leq C$ for all $u\in T_{\xo}X$, $\|u\|\leq 1$. Therefore,
\begin{equation}
C' \leq L_{\xo}^{(X)}(u) \leq C   \:\:\:\:\;\:, \:\:\: u\in T_{\xo}X,\; \|u\|\leq 1\;,
\end{equation}
for some other global constant $C' > 0$ (depending on $K$ but not an $\lambda$ and $u$).

\def\naf{\nabla f(z)}
\def\so{\sigma_0}
 \def\tp{\tilde{p}}
 
Consider the map $\Phi : U \longrightarrow S$ near $\sigma_0 = (\xo,\xio)$ defined by
$\Phi(x,\xi)  =  {\iii} \circ \phi_{2\lambda}(x,\xi)$.
In fact by the same formula (see below for more details) one defines $\Phi$ as a local
smooth  map $\Phi: T(\R^n) = \R^n \times \R^n \longrightarrow T(\R^n)$ near $\sigma_0$. 
Given $\ep > 0$, we will assume $\delta \in (0,\ep]$ is chosen sufficiently small , so that
$\Phi$ is well-defined and $\Phi (U) \subset S$.  
Moreover,  $\Phi(z) = z'$ implies $\Phi(z') = z$ (whenever $\Phi(z')$ is defined) and 
locally $\Phi(W^u_\ep(z)) = W^s_\ep(z')$. Finally, it is important to remark that $\Phi$ preserves the set $\mt$. 
Indeed, $\phi_{2\lambda} (\mt) = \mt$ and $\iii(\mt) = \mt$, as well. So, in particular
\be
\Phi(U \cap \mt) \subset  S\cap \mt\;.
\ee

To write down a more explicit expression for $\Phi$, let
$f$ be a defining function for $\dk$ in a neighbourhood of $\dk$ in $\R^n$  so that
$\| \nabla f \| = 1$ near $\zo$ and $\nabla f(z) = \nu(z)$ is the outward unit normal to $\dk$ at $z\in \dk$.
Then $\dk = f^{-1}(0)$ (locally near $\zo$). Given $(x,\xi)\in \R^n\times\R^n$ close to $\sigma_0$, there
exist a unique $z(x,\xi)\in \dk$ and a unique minimal $t(x,\xi)\in \R_+$ with 
$z(x,\xi) = x + t(x,\xi)\xi\in \dk\;,$
i.e. such that 
\begin{equation}
f(x + t(x,\xi)\xi) = 0\;.
\end{equation}
By $\eta(x,\xi)$ we denote the reflection of $\xi$ with respect to $\nabla f(z(x,\xi))$, i.e.
$$\eta(x,\xi) = \xi - 2 \la \xi , \nabla f(z)\ra \,\naf\;.$$
Here and in what follows we denote for brevity $z = z(x,\xi)$. We will also use the notation $t = t(x,\xi)$
and  $\eta = \eta(x,\xi)$. We then have 
$$\Phi(x,\xi) = (g(x,\xi), -\eta(x,\xi))\;,$$
where
$g(x,\xi) = z + (2\lambda - t)\eta \;.$
Since $\Phi(U) \subset S$  and $\Phi$ is a local diffeomorphism between $U$ and $S$,
we have $d\Phi_{\sigma} (E^u(\sigma)) =  E^s(\sigma)$ for every $\sigma \in U$. 
Moreover, it is easy to see that $d\Phi_\sigma$ preserves the sets of tangent vectors
to $\mt$, namely if $\sigma \in U\cap \mt$ and $\xi \in E^u(\sigma)\setminus \{0\}$ is tangent to $\mt$ at 
$\sigma$, then  $d\Phi_\sigma \cdot \xi$  is tangent to $\mt$ at $\Phi(\sigma)$.

It is well known that we can take the  constant $C> 0$ so large that $\|d\Phi_{\sigma}\| \leq C$
for any $\sigma\in U$ and any choice of $\so$ (see e.g. \cite{kn:Si2}, \cite{kn:Ch1} or \cite{kn:BCST}).
(See also the proof of Lemma 3.1 in Sect. 4 for an explicit formula for $d \Phi_{\so}$.)

Set $L = L_{\xo}^{(X)}$ and $H = L_{\zo}$, and for any $u \in T_{\xo}X$,  consider the vectors
$$v(u) = (u, L u)\in E^u(\sigma_0)\quad, \quad w(u) = d\Phi (\sigma_0)\cdot v(u)\in E^s(\sigma_0)\;.$$
It is easy to see that
\be
\|v(u) \| \leq \sqrt{1+ C^2}\, \|u\| \;.
\ee

The following lemma is the main technical ingredient in the proofs of Theorem 1.1 and 
Proposition 6.2. Its proof is given in Sect. 4 below.
 
\bigskip

\noindent
{\sc Lemma 3.1.}  {\it For any $u,u'\in T_{\xo}X$ we have
$d\alpha (v(u), w(u')) =  \left\langle u, P\,u'\right\rangle \;,$
where the linear operator $P$ is given by 
$P = 2H + 2L + 2\lambda (HL + LH + L^2 + \lambda \, L H L) \;.$
Consequently, if  $\kappa > 0$ is the minimal principal curvature at a point on $\dk$ and $\ep$ and $\lambda$ are chosen sufficiently small, 
then $P$ is positive definite, $\la u, Pu\ra \geq \kappa\, \|u\|^2$, and therefore
$|d\alpha (v(u), w(u))| \geq \kappa\, \|u\|^2$ for all $u \in T_{\xo}X$. }

\bs

\noindent
{\it Proof of Theorem 1.1.} Take $\ep > 0$ sufficiently small and then $\lambda$ with (3.1) 
small enough so that
the operator $P$ in Lemma 3.1 is positive definite, where $\so$ is chosen as above.
(Notice that $H$ and $L$ are uniformly bounded from below and above regardless of the 
choice of the point $\so$.)
More precisely, as stated in Lemma 3.1, if $\kappa > 0$ is the minimal principal 
curvature at a point on $\dk$, we can choose 
$0 < \ep < \lambda$ so small that $\la u,Pu\ra \geq \kappa\, \|u\|^2$ for any $u\in T_{\xo}X$. 

Set $z_0 = \so$ and let $b\in \hE^u_\mt(z_0)$, $\|b\| = 1$. Then $b = v(u)$ for some 
$u \in T_{\xo}X$. Moreover, by (3.3),
$w(u) = d\Phi(\so)\cdot v(u) \in \hE^s_\mt(z_0)$, and so $a = w(u)/\|w(u)\|$ is a 
unit vector in $\hE^s_\mt(z_0)$.
By Lemma 3.1,
$$|d\alpha(a,b)| = \frac{1}{\|w(u)\|}\, |d\alpha(v(u), w(u))| \geq \frac{\kappa\, \|u\|^2}{\|w(u)\|}
\geq \frac{\kappa\, \|u\|^2}{C\, \|v(u)\|} =  \frac{\kappa\, \|u\|^2}{C}\;.$$
On the other hand (3.5) implies
$1 = \|b\| = \|v(u)\| \leq \sqrt{1+C^2}\, \|u\|$, so $|d\alpha(a,b)| \geq \frac{\kappa}{C\, (1+C^2)}$.
\endofproof

\bs

As a consequence of Lemma 3.1  one can also derive the following which however we do not need in this paper.

\bs

\noindent
{\sc Proposition 3.2.} {\it  For every $z\in \mt$ and every $\delta > 0$ there exists $\tz \in \mt\cap W^u_\delta(z)$ such that for any
non-zero tangent vector $b\in E^u(\tz)$ to $\mt$ there exists a tangent vector $a\in E^s(\tz)$ to $\mt$ with $d\alpha (a,b) \neq 0$}.

\section{Proof of Lemma 3.1}
\setcounter{equation}{0}

 We will use the notation from Sect. 3.
Recall that the standard symplectic form $d\alpha$ has the form
$$d\alpha(\,(u,\tilde{u})\, , \, (p,\tilde{p})\,) = - \la u, \tilde{p}\ra + \la \tilde{u} ,p\ra\;,$$
where $(u,\tilde{u}), (p,\tilde{p}) \in T^*(\R^n)$.
Given $u,u' \in T_{\xo}X$, let  $v(u) = (u,\tu)\in E^u(\sigma)$,
$v(u') = (u',\tu')\in E^u(\sigma)$ and $w(u') = d\Phi_{\sigma_0}(v(u')) = (p,\tp)$.
Then
$$\pmatrix{p\cr \tp\cr} = \pmatrix{\partial_x g(\sigma) & \partial_\xi g(\sigma)\cr
-\partial_x\eta (\sigma) & -\partial_\xi \eta(\sigma)\cr}\;\pmatrix{u'\cr \tu'\cr}
= \pmatrix{ \partial_x g(\sigma)\, u' + \partial_\xi g(\sigma)\, \tu'\cr
-\partial_x\eta (\sigma) \,u' - \partial_\xi \eta(\sigma)\,\tu'\cr}\;,$$
and so
\begin{eqnarray}
d\alpha(v(u), w(u'))  
& = &    \la \,u\,, \,\partial_x \eta(\sigma)\, u'\, \ra  
+ \la\, u\,,\,\partial_\xi \eta(\sigma)\, \tu'\,\ra\\
&    &       \:\:\: + \la \,\tu \,,\, \partial_x g (\sigma) \,u'\,\ra 
 + \la \,\tu\,,\, \partial_\xi g(\sigma)\,\tu'\, \ra\nonumber\;. 
\end{eqnarray}

One needs the derivatives of $g$ and $\eta$. Differentiating (3.4) gives
\begin{equation}
\nabla_x t(\sigma_0) = \nabla f(\zo) = - \xio \:\:\:, \:\:\:\:
\nabla_\xi t(\sigma_0) = \lambda \,\nabla f(\zo) = - \lambda\, \xio\;.
\end{equation}
Moreover, $z = x + t\xi$ implies
\begin{equation}
\frac{\partial z_\ell}{\partial x_j} (x,\xi) = \delta_{j\ell} + \frac{\partial t}{\partial x_j}(x,\xi)
\,\xi_\ell \:\:\:, \:\:\:
\frac{\partial z_\ell}{\partial \xi_j} (x,\xi) = t \delta_{j\ell} + \frac{\partial t}{\partial \xi_j}(x,\xi) \,\xi_\ell\;.
\end{equation}

Next, we have
$g(x,\xi) =  z+ (2\lambda - t)\, \eta = x + 2\lambda \, \xi - 2(2\lambda - t)\, \la \xi , \nabla f(z)\ra \,\naf\;.$
Hence
\begin{eqnarray*}
    \frac{\partial g_i}{\partial x_j}(x,\xi) 
& = & \delta_{ij}  + 2 \frac{\partial t}{\partial x_j} (x,\xi)\, \la \xi,\naf\ra \, \frac{\partial f}{\partial x_i}(z)
     - 2(2\lambda - t)\left(\sum_{k,\ell=1}^n \xi_k \frac{\partial^2 f}{\partial x_k\partial x_\ell}(z)\;
\frac{\partial z_\ell}{\partial x_j}\right)\, \frac{\partial f}{\partial x_i}(z)\\
&   & - 2(2\lambda - t)\la \xi , \naf\ra\,  \sum_{\ell=1}^n \frac{\partial^2 f}{\partial x_i\partial x_\ell}(z)\;
\frac{\partial z_\ell}{\partial x_j}\;,
\end{eqnarray*}
\begin{eqnarray*}
&&\frac{\partial g_i}{\partial \xi_j}(x,\xi) 
 =  2\lambda  \, \delta_{ij} + 2 \frac{\partial t}{\partial \xi_j} (x,\xi) \,\la \xi,\naf\ra 
\frac{\partial f}{\partial x_i}(z) - 2(2\lambda -t)\, \frac{\partial f}{\partial x_j}(z)
\frac{\partial f}{\partial x_i}(z) \\
& & - 2(2\lambda-t)\left(\sum_{k,\ell=1}^n \xi_k  \frac{\partial^2 f}{\partial x_k\partial x_\ell}(z)\;
\frac{\partial z_\ell}{\partial \xi_j}\right)\, \frac{\partial f}{\partial x_i}(z)
    - 2(2\lambda -t)\la \xi , \naf\ra\, 
\sum_{\ell=1}^n \frac{\partial^2 f}{\partial x_i\partial x_\ell}(z)\; \frac{\partial z_\ell}{\partial \xi_j}\;,
\end{eqnarray*}
Similarly,
\beqn
\frac{\partial \eta_i}{\partial x_j}(x,\xi) 
& = &  - 2\left(\sum_{k,\ell=1}^n \xi_k \frac{\partial^2 f}{\partial x_k\partial x_\ell}(z)\;
\frac{\partial z_\ell}{\partial x_j}\right)\, \frac{\partial f}{\partial x_i}(z)
   - 2\la \xi , \naf\ra\,  \sum_{\ell=1}^n \frac{\partial^2 f}{\partial x_i\partial x_\ell}(z)\; \frac{\partial z_\ell}{\partial x_j}\;,
\eeqn
\beqn
\frac{\partial\eta_i}{\partial \xi_j}(x,\xi) 
& = &  \delta_{ij} - 2\, \frac{\partial f}{\partial x_j}(z) \frac{\partial f}{\partial x_i}(z) - 2\, \left(\sum_{k,\ell=1}^n \xi_k 
\frac{\partial^2 f}{\partial x_k\partial x_\ell}(z)\; \frac{\partial z_\ell}{\partial \xi_j}\right)\, \frac{\partial f}{\partial x_i}(z)\\
&   &  - 2\,\la \xi , \naf\ra\, \sum_{\ell=1}^n \frac{\partial^2 f}{\partial x_i\partial x_\ell}(z)\;
\frac{\partial z_\ell}{\partial \xi_j}\;.
\eeqn

Notice that for $v = (u,\tu)\in E^u(\sigma_0)$ we have $u, \tu \perp \xi_0$. Therefore (4.2) imply
$\la \,\nabla_xt(\sigma_0)\,, \,u\,\ra = \la \,\nabla_\xi t(\sigma_0)\,,\,  u\,\ra = 0$,
and the same holds with $u$ replaced by $\tu$. This and (4.3) give
$$\sum_{j=1}^n \frac{\partial z_\ell}{\partial x_j}(\sigma_0)\, u_j = u_\ell\quad , \quad 
\sum_{j=1}^n \frac{\partial z_\ell}{\partial \xi_j}(\sigma_0)\, u_j = \lambda\, u_\ell$$
for any $u \in T_{\xo}X$.  Using these, $t(\sigma_0) = \lambda$, the above formulae and the Hessian matrix\\
$\di H' = \left( \frac{\partial^2 f}{\partial x_k \partial x_\ell} (\zo)\right)_{k,\ell=1}^n\;,$ 
one gets :
\beqn
\sum_{j=1}^n \frac{\partial g_i}{\partial x_j}(\sigma_0) \, u'_j
& = &  u'_i - 2\lambda \left(\sum_{k,\ell=1}^n \xio_k u'_\ell \, \frac{\partial^2 f}{\partial x_k \partial
x_\ell} (\zo)\right)\, \frac{\partial f}{\partial x_i}(\zo) + 2\lambda\,\sum_{\ell=1}^n u'_\ell \,
\frac{\partial^2 f}{\partial x_i \partial x_\ell} (\zo)\\
& = & u'_i - 2\lambda \la\, u'\,,\,H'\xio\,\ra \, \frac{\partial f}{\partial x_i}(\zo) + 2\lambda\, (H'u')_i
= u'_i + 2\lambda\, (H'u')_i\;,
\eeqn 
where $(H'u')_i$ is the $i$th coordinate of the (column) vector $H'u'$. Here we used the fact that
$\xio = - \nabla f(\zo)$ and $H'\nabla f(\zo) = 0$, since $\| \nabla f\| = 1$ near $\dk$. Similarly,
\beqn
\sum_{j=1}^n \frac{\partial g_i}{\partial \xi_j}(\sigma_0) \, \tu'_j
& = &  2\lambda \tu'_i - 2\lambda^2\left(\sum_{k,\ell=1}^n \xio_k \tu'_\ell \, \frac{\partial^2 f}{\partial x_k
\partial x_\ell} (\zo)\right)\, \frac{\partial f}{\partial x_i}(\zo)\\ 
&   & + 2\lambda^2\,\sum_{\ell=1}^n \tu'_\ell \, \frac{\partial^2 f}{\partial x_i \partial x_\ell} (\zo)
 = 2\lambda \tu'_i + 2\lambda^2 (H' \tu')_i\;,
\eeqn 
\beqn
\sum_{j=1}^n \frac{\partial \eta_i}{\partial x_j}(\sigma_0) \, u'_j
 =   - 2 \left(\sum_{k,\ell=1}^n \xio_k u'_\ell \, \frac{\partial^2 f}{\partial x_k \partial x_\ell}
(\zo)\right)\, \frac{\partial f}{\partial x_i}(\zo) 
+ 2\,\sum_{\ell=1}^n u'_\ell \, \frac{\partial^2 f}{\partial x_i \partial x_\ell} (\zo) = 2(H'u')_i\;,
\eeqn 
\beqn
\sum_{j=1}^n \frac{\partial \eta_i}{\partial \xi_j}(\sigma_0) \, \tu'_j
& = & \tu'_i - 2\lambda \left(\sum_{k,\ell=1}^n \xio_k \tu'_\ell \, \frac{\partial^2 f}{\partial x_k \partial x_\ell} (\zo)\right)\, 
\frac{\partial f}{\partial x_i}(\zo) + 2\lambda\,\sum_{\ell=1}^n \tu'_\ell \, \frac{\partial^2 f}{\partial x_i \partial x_\ell} (\zo)\\
& = & \tu'_i + 2\lambda\, (H'\tu')_i\;.
\eeqn 
The last four formulae imply
$\la \,u\,, \,\partial_x \eta(\sigma_0)\, u'\, \ra = 
\sum_{i,j=1}^n u_iu'_j\, \frac{\partial \eta_i}{\partial x_j}(\sigma_0) = 2\la\, u\,,\, H'u'\,\ra\;,$
and similarly
$\la\, u\,,\,\partial_\xi \eta(\sigma_0)\, \tu'\,\ra  =
\la u,\tu'\ra + 2\lambda \la\, u\,,\, H'\tu'\, \ra$,
$\la \,\tu \,,\, \partial_x g (\sigma_0) \,u'\,\ra =
\la \tu , u'\ra + 2\lambda \la \tu , H'u'\ra$, and
$ \la \,\tu\,,\, \partial_\xi g(\sigma_0)\,\tu'\, \ra =
2 \lambda \la \tu , \tu'\ra + 2\lambda^2 \la \tu , H'\tu'\ra$.
Combining these with (4.1), one gets
$$d\alpha(v(u), w(u')) = 2\la\, u\,,\, H'u'\,\ra + \la u, \tu'\ra + 2\lambda \la u, H'\tu'\ra +
\la \tu,u'\ra + 2\lambda \la\, \tu\,,\, H'u'\, \ra + 2 \lambda \la \tu , \tu'\ra + 2\lambda^2 \la \tu , H'\tu'\ra\;.$$
Using $\tu = L(u)$, $\tu' = L(u')$ and the fact that $H'u = Hu$ for all $u \in T_{\xo}X$, it now follows that
\begin{eqnarray*}
d\alpha(v(u), w(u')) 
& = &  2\la\, u\,,\, Hu'\,\ra + \la u, L u'\ra + 2\lambda\la u, HL u'\ra + \la L u,u'\ra + 2\lambda \la\, L u\,,\, Hu'\, \ra \\
&    & \:\: + 2 \lambda \la L u , L u'\ra + 2\lambda^2 \la L u , H L u'\ra =  \la u, Pu'\ra\;,
\end{eqnarray*}
where $P = 2H + 2L + 2\lambda \, (HL + LH + L^2 + \lambda\, LHL)$.
\endofproof

\section{Pinched open billiard flows}
\setcounter{equation}{0}

In this section we describe some open billiard flows in $\R^n$ ($n\geq 3$) that satisfy the pinching
condition (P). (Clearly open billiards in $\R^2$ always satisfy this condition.) As one can see below, the
estimates we use are rather crude, so one would expect that with more sophisticated methods larger
classes of open billiard flows could be shown to satisfy the condition (P).

First, we derive a formula which is useful in getting estimates for $\|d\phi_t(x)\cdot u\|$ ($u \in E^u(x)$, $x\in \mt$),
both from above and below.
From the arguments in this section one can also derive a representation for the Jacobi fields along a 
billiard trajectory.

In what follows we use the notation from the beginning of Sect. 3. Here {\bf we assume that the boundary $\dk$ is
at least $C^3$ smooth}.

Fix  for a moment a point $x_0 = (q_0,v_0) \in \mt_\ep$. If $\ep > 0$ is sufficiently small, then 
$W^u_\ep(x_0)$ has the form (cf. \cite{kn:Si1},\cite{kn:Si2}) 
$W^u_\ep(x_0) = \{ (x,\nu_X(x)) : x\in X\}$
for some smooth  hypersurface $X$ in $\R^n$ containing the point $q_0$ such that $X$
is strictly convex with respect to the unit normal field $\nu_X$. Denote by
$B_x : T_q X \longrightarrow T_qX$ the {\it curvature operator} ({\it second fundamental form})
of $X$ at $q \in X$. Then $B_x$ is positive definite with respect to the normal field $\nu_X$ (\cite{kn:Si2}).

Given a point $q\in X$, let $\gamma(x)$ be the {\it forward billiard
trajectory}  generated by $x = (q,\nu_X(q))$. Let $q_1(x), q_2(x), \ldots$ 
be the reflection points of this trajectory 
and let $\xi_j(x)\in \sn$ be the {\it reflected direction} of $\gamma(x)$ at $q_j(x)$.
Set $q_0(x) = q$, $t_0(x) = 0$ and denote by $t_1(x), t_2(x), \ldots$ 
the times of the consecutive reflections of the trajectory $\gamma(x)$ at $\dk$. Then 
$t_j(x) = d_0(x) + d_1(x) + \ldots + d_{j-1}(x)$, where
$d_j(x) = \| q_{j+1}(x) - q_j(x)\| \:\:\: , \:\:\: 1\leq j\;.$
Given $t \geq 0$, denote by $u_t(q)$ the {\it shift} of $q$ along the
trajectory $\gamma(x)$ after time $t$. Set
$X_t = \{ u_t(q) : q\in X\}\;.$
When $u_t(q)$ is not a reflection point of $\gamma(x)$, then locally
near $u_t(q)$, $X_t$ is a smooth convex $(n-1)$-dimensional surface in $\R^n$ with "outward" 
unit normal given by the direction $v_t(q)$ of $\gamma(x)$ at $u_t(q)$ (cf. \cite{kn:Si2}).

Fix for a moment $t > 0$ such that $t_m(x_0) < t < t_{m+1}(x_0)$ for some $m \geq 1$,
and assume that $q(s)$, $0 \leq s\leq a$, is a $C^3$ curve on $X$ with $q(0) = q_0$ 
such that for every $s \in [0,a]$ we have $t_m(x(s) ) < t < t_{m+1}(x(s))$,
where $x(s) = (q(s), \nu_X(q(s)))$.  Assume also that $a > 0$ is so small that for all
 $j = 1,2,\ldots,m$ the reflection points $q_j(s) = q_j(x(s))$ belong to
the same boundary component $\dk_{i_j}$ for every $s\in [0,a]$.

We will now estimate $\| d\phi_t(x_0)\cdot \xi_0\|$, where $\xi_0 = \dot{q}(0) \in T_{q_0}X$.

Clearly $\phi_t(x(s)) = (p(s), v_m(x(s)))$, where $p(s)$, $0\leq s\leq a$, is a $C^3$ curve on $X_t$.
For brevity denote by $\gamma(s)$ the {\it forward billiard trajectory}  
generated by $(q(s), \nu(q(s)))$ and set $q_0(s) = q(s)$.
Let $\xi_j(s)\in S^{n-1}$ be the {\it reflected direction} of $\gamma(s)$ at $q_j(s)$
and let $\varphi_j(s)$ be the {\it angle} between $\xi_j(s)$ and 
the outward unit normal $\nu(q_j(s))$ of $\dk$ at $q_j(s)$.
Let $\phi_t(x(s)) = (u_t(s), v_t(s))$, and let $ t_j(s) = t_j(x(s))$
be the times of the consecutive reflections of the trajectory $\gamma(s)$ at $\dk$. 
Set $d_j(s) = d_j(x(s)) =  \| q_{j+1}(s) - q_j(s)\|$
($0\leq j\leq m-1$),  $t_0(s) = 0$,  $t_{m+1}(s) = t$ and $d_{m}(s) = t - t_m(s)$.  
Denote by $k_t(s)$ the {\it normal curvature} of $X_t$ at $u_t(s)$ in the direction of $\frac{d}{ds}u_t(s)$.

Next, let $k_0(s)$ be the {\it normal curvature} of $X$ at $q(s)$ in the direction
of $\dot{q}(s)$, and for $j > 0$ let $k_j(s) > 0$ be the {\it normal curvature} of 
$X_{t_j(s)} = \lim_{t\searrow t_j(s)} X_t$ at $q_j(s)$ in the direction $\hat{u}_j(s)$ ($\| \hat{u}_j(s)\| = 1$)
of $\lim_{t\searrow t_j(s)}\frac{d}{ds'}\left(u_{t}(s')\right)_{|s' = s}$. For $j \geq 0$ let
$$B_j(s): T_{q_j(s)}(X_{t_j(s)}) \longrightarrow T_{q_j(s)}(X_{t_j(s)})$$
be the {\it curvature operator} (second fundamental form) of $X_{t_j(s)}$ at $q_j(s)$, and define $\ell_j(s) > 0$ by
\be
[1+ d_j(s)\ell_j(s)]^2 = 1+ 2d_j(s)k_j(s) + (d_j(s))^2\, \|B_j(s)\hat{u}_j(s)\|^2\;.
\ee
Finally, set
\begin{equation}
\delta_j(s) = \frac{1}{ 1 + d_j(s) \ell_j(s)} \:\:\:\: , \:\:\: 0\leq j \leq m\;.
\end{equation}

\ms

\noindent
{\sc Proposition 5.1.} {\it For all $s\in [0,a]$  we have}
\begin{equation}
\| \dot{q}_0(s)\| = \| \dot{p}(s)\| \delta_{0}(s)\delta_{1}(s) \ldots \delta_m(s)\;. 
\end{equation}

\ms

As we mentioned in the Introduction, the above formula can be easily derived from the more general 
study of the evolution of unstable fronts in multidimensional dispersing billiards in \cite{kn:BCST} (see Section 5 there).
Apart from that, one could prove (5.3) by using a simple modification of the argument in the Appendix of \cite{kn:St2}
dealing with the two-dimensional case. We omit the details. 

We will now use Proposition 5.1 to prove Proposition 1.2.

In the notation above, let  $q_j = q_j(x)$ be the reflection points of the billiard trajectory $\gamma(x)$
for some $x = (q,\nu_X(q))$, with $q\in X$, and let $t_j = t_j(x)$ and $d_j = d_j(x)$. Consider 
the curvature operator
$B_j = B_{q_j} : \Pi_j = T_{q_j}(X_{t_j}) \longrightarrow \Pi_j\;,$
and let $S_j: \R^n \longrightarrow \R^n$  be the {\it symmetry} with 
respect to the {\it tangent space} $\tt_j = T_{q_j}(\dk)$; notice that
$S_j(\Pi_{j-1}) = \Pi_j$. 
Let $N_j : \tt_j \longrightarrow \tt_j$  be the {\it curvature operator} (second fundamental form) of $\dk$ at $q_j$.

Notice that $\Pi_j$ is the hyperplane in $\R^n$ passing through $q_j$
and orthogonal to $\xi_j = v_{t_j(x)}(x)$; it will be identified with the $(n-1)$-dimensional
vector subspace of $\R^n$ orthogonal to $\xi_j$.  

Before going on we need to recall the representation of the operator $B_j$
due to Sinai \cite{kn:Si2} (cf. also Chernov \cite{kn:Ch1}).  Introduce the linear maps 
$V_j: \Pi_j \longrightarrow \tt_j $, $V^*_j : \tt_j \longrightarrow \Pi_j$
where $V_j$ is (the restriction to $\Pi_j$ of) the projection to $\tt_j$ along the
vector $\xi_j$, while $V^*_j$ is the projection to $\Pi_j$ along the normal vector
$\nu_j = \nu (q_j)$. (Considering $V_j : \R^n \longrightarrow \tt_j$ and
$V^*_j : \R^n \longrightarrow \Pi_j$, $V^*_j$ is the self-adjoint of $V_j$.) Let $\varphi_j$ be the {\it angle} between 
$\nu_j$ and $\xi_j$. Then (\cite{kn:Si2})
\begin{equation}
B_j = S_j\, B^-_{j}\,S_j + 2\cos\varphi_j\, V^*_j N_j V_j
\end{equation}
for $1\leq j \leq m$, 
$B^-_j = B_{j-1} \, (I + d_{j-1}\, B_{j-1})^{-1}\;,$
and $B_{m+1} = B_{m} \, (I + t'\, B_{m})^{-1}\;,$ where $t' = t - t_m \geq \ep$.

Let $\mu_j(x_0) \leq \lambda_j(x_0)$ be the minimal and the maximal eigenvalues of the operator $B_j$.
If  $\lambda$ is an eigenvalue of $ B_{j-1}$, then 
$\lambda/(1+d_{j-1}\lambda)$ is an eigenvalue of $B^-_j $, and
$\frac{\lambda}{1+ \lambda \, d_{j-1}} = \frac{1}{1/\lambda + d_{j-1}} < \frac{1}{d_0}\;,$
Next, a simple calculation shows that the spectrum of the operator $V^*_j N_j V_j$ lies
in the interval $[ \kmin, \frac{\kmax}{\cos^2\varphi_j}]$.
Thus, using  (5.4) we get
\be
\mu_0 \leq  2 \cos\varphi_j \, \kmin \leq \mu_j(x_0) \leq \lambda_j(x_0) \leq \frac{1}{d_0} + \frac{ 2\kmax}{\cos\varphi_j}  \leq \lambda_0\;,
\ee
where $\mu_0$ and $\lambda_0$ are as in Sect. 1.

\bs

\noindent
{\it Proof of Proposition} 1.2. 
Before we continue, notice that there exist global constants $0 < c_1 < c_2$ such that
$c_1\, \|\xi\| \leq \|u\| \leq c_2\, \|\xi\|$ for any $u = (\xi, \eta) \in E^u(x)$, $x\in \mt$ 
(see formula (3.5) above).

Assume that (1.1) holds.
Fix an arbitrary $x_0 = (q_0,v_0) \in \mt_\ep$ and $t >0$. We will now use the notation from the
beginning of this section. 

To estimate $\| d\phi_t(x_0)\cdot u\|$ for a given unit vector $u = (\xi, \eta) \in E^u(x_0)$, 
consider a $C^1$ curve $q(s)$, $0 \leq s\leq a$, on $X = \pr_1(W^u_\ep(x_0))$ with 
$q(0) = q_0$ and $\dot{q}(0) =\xi$, and define $q_j(s)$, $j= 1, \ldots,m$ and $p(s)$ as in the
beginning of this section. Then $p(s) = \pr_1(d\phi_t(x(s)))$, so
$c_1\|\dot{p}(0)\| \leq \|d\phi_t(x_0)\cdot u\| \leq c_2\|\dot{p}(0)\|$. Using this and Proposition 5.1, we get 
\be
\frac{c_1\, \|u\|}{ c_2\, \delta_{1}(0)\delta_{2}(0) \ldots \delta_m(0)} \leq \|d\phi_t(x_0)\cdot u\| \leq 
\frac{c_2\, \|u\|}{ c_1\, \delta_{1}(0)\delta_{2}(0) \ldots \delta_m(0)}\;.
\ee

Recall that each $\delta_j$ is given by (5.2) and (5.1), so if $0 < \mu_j(x_0) \leq \lambda_j(x_0)$ are the minimal and
maximal eigenvalues of the operator $B_j(0)$, then 
$$(1+ d_j(0)\ell_j(0))^2 \leq 1 + 2d_j(0)\lambda_j + d_j^2(0) \, \lambda_j^2 = (1+d_j(0)\, \lambda_j(x_0))^2\;,$$
so $\ell_j(0) \leq \lambda_j(x_0)$. Similarly, $\mu_j(x_0) \leq \ell_j(0)$. Moreover, it follows from (5.5)
that $\mu_0 \leq \mu_j(x_0)$ and $\lambda_j(x_0) \leq \lambda_0$ for all $j\geq 1$ and all $x_0\in \mt$.

Assuming $\|u\| = 1$ and recalling that $d_j(x_0) = d_j(0)$ and $t > d_1(x_0) + \ldots + d_m(x_0)$, 
(5.6) and (5.2)  give
\begin{eqnarray*}
\frac{1}{t} \ln \|d\phi_t(x_0)\cdot u\| 
& \leq & \frac{\ln (c_2/c_1)}{t} + \frac{1}{t} \sum_{j=1}^m \ln (1+ d_j(x_0)\, \lambda_j(x_0))\\
& \leq & \frac{\ln (c_2/c_1)}{t} + \frac{\sum_{j=1}^m \ln (1+ (d_0+a)\, \lambda_0)}{d_1(x_0) 
+ \ldots + d_m(x_0)}\\
& \leq & \frac{\ln (c_2/c_1)}{t} + \frac{m \ln (1+ (d_0+a)\, \lambda_0)}{m d_0}\\ 
& \leq & \frac{\ln (c_2/c_1)}{t} +  \frac{ \ln (1+ (d_0+a)\, \lambda_0)}{d_0} =  
\frac{\ln (c_2/c_1)}{t} + \beta_0\;,\\
\end{eqnarray*}
so $\|d\phi_t(x_0)\cdot u\| \leq (c_2/c_1)\, e^{t\,\beta_0}$ for all  $ t > 0$.

In a similar way, using (5.6) one derives that 
$\|d\phi_t(x_0)\cdot u\| \geq c'\, (c_1/c_2)\, e^{t\,\alpha_0}$ for  $t > 0$, 
where $\alpha_0 =  \frac{ \ln (1+ d_0\, \mu_0)}{d_0 +a}$ and $c' > 0$ is another global constant. 
Finally, notice that (1.1) implies $2 \alpha_0 \geq \beta_0 + \alpha$ for some global constant $\alpha > 0$.
Hence the condition (P) is satisfied.
\endofproof

\def\hhu{\hh^{u}}
\def\hhs{\hh^{s}}
\def\Lu{L^{u}}
\def\Ls{L^{s}}
\def\ta{\tilde{a}}

\section{Dolgopyat type estimates for pinched open billiard flows}
\setcounter{equation}{0}

Let $\phi_t : M \longrightarrow M$ be a  $C^1$ flow on complete (not necessarily compact) Riemann manifold  $M$, and let
$\mt$ be a basic set for $\phi_t$. 
It follows from the hyperbolicity of $\mt$  that if  $\epsilon > 0$ is sufficiently small, there exists $\delta > 0$ such that if 
$x,y\in \mt$ and $d (x,y) < \delta$, then $\wloc^s(x)$ and $\phi_{[-\ep,\ep]}(\wloc^u(y))$ intersect at exactly 
one point $[x,y ] \in \mt$  (cf. \cite{kn:KH}). That is, there exists a unique $t\in [-\ep, \ep]$ such that $\phi_t([x,y]) \in \wloc^u(y)$. 

Let $\rr = \{R_i\}_{i=1}^k$ be a Markov family for $\phi_t$ over $\mt$ consisting of 
rectangles $R_i = [U_i ,S_i ]$, where $U_i$ (resp. $S_i$) are (admissible) subsets of $W^u_{\ep}(z_i) \cap \mt$
(resp. $W^s_{\ep}(z_i) \cap \mt$) for some $\ep > 0$ and $z_i\in \mt$ (cf. e.g. \cite{kn:PP}  for details; see also \cite{kn:D}). 
The {\it first return time function}  $\tau : R = \cup_{i=1}^k R_i \longrightarrow [0,\infty)$ and the standard Poincar\'e map 
$\pp: R \longrightarrow R$ are 
Lipschitz when restricted to an appropriate large subset of $R$.
Set $U = \cup_{i=1}^k U_i$ and define the {\it shift map} $\sigma : U \longrightarrow U$ 
by $\sigma = p\circ \pp$, where $p : R \longrightarrow U$
is the projection along the leaves of local stable manifolds. 
Let $\hU$ be the set of all $u \in U$ whose orbits do not 
have common points with the boundary of $R$ (in $\mt$). Given a Lispchitz function (or map) 
on $\hU$, we will identify it with its (unique) Lipschitz extension to $U$.
Assuming that the local stable and unstable laminations over $\mt$ are Lipschitz, 
the map $\sigma$ is {\it essentially Lipschitz} on $U$ 
in the sense that there exists a constant $L > 0$ such that if $x,y \in U_i \cap \sigma^{-1}(U_j)$ 
for some $i,j$, then $d(\sigma(x), \sigma(y)) \leq L\, d(x,y)$.
The same applies to $\tau : U \longrightarrow \R$.

Given a Lipschitz real-valued function $f$  on $\hU$, set $g = g_f = f - P\tau$, where $P = P_f\in \R$ is the unique 
number such that the topological pressure $\Pr_\sigma(g)$ of $g$ with respect to $\sigma$ is zero (cf. e.g. \cite{kn:PP}). 
For $a, b\in \R$, one defines the {\it Ruelle operator} $L_{g-(a+\i b)\tau} : \clip (\hU) \longrightarrow \clip (\hU)$ in 
the usual way (cf. e.g. \cite{kn:PP} or \cite{kn:D}),
where $\clip (\hU)$ is the space of Lipschitz functions $g: \hU \longrightarrow \C$. By  $\Lip(g)$ we denote the Lipschitz constant of $g$
and  by $\| g\|_0$ the {\it standard $\sup$ norm}  of $g$ on $\hU$.

We will say  that the {\it Ruelle transfer operators related to the function $f$ on $U$ are eventually contracting} 
if for every $\epsilon > 0$ there exist constants $0 < \rho < 1$, $a_0 > 0$ and  $C > 0$ such that  if $a,b\in \R$  are such that 
$|a| \leq a_0$ and $|b| \geq 1/a_0$, then for every integer $m > 0$ and every  $h\in \clip (\hU)$ we have
$$\|L_{f -(P_f+a+ \i b)\tau}^m h \|_{\lip,b} \leq C \;\rho^m \;|b|^{\ep}\; \| h\|_{\lip,b}\; ,$$
where the norm $\|.\|_{\lip,b}$ on $\clip (\hU)$ is defined by $\| h\|_{\lip,b} = \|h\|_0 + \frac{\lip(h)}{|b|}$.
This implies in particular that the spectral radius  of $L_{f-(P_f+ a+\i b)\tau}$ in $\clip(\hU)$ does not exceed  $\rho$.

Next, assume that $\phi_t$ is a $C^2$ contact flow on $M$ with a $C^2$ invariant  contact
form $\omega$. The following condition says  that $d\omega$ is in some sense 
non-degenerate on $\mt$ near some of its points:

\bs

\def\ty{\tilde{y}}

\noindent
{\sc (ND)}:  {\it There exist $z_0\in \mt$, $\delta_0 > 0$  and  $\mu_0 > 0$ such that
for any $\delta \in (0,\delta_0]$, any $\hz \in \mt \cap W^u_{\delta}(z_0)$ and any unit vector 
$b \in E^u(\hz)$ tangent to 
$\mt$ at $\hz$ there exist $\tz \in \mt \cap W^u_{\delta}(\hz)$, $\ty \in W^s_\delta(\tz)$ and 
a unit vector $a \in E^s(\ty)$ tangent to $\mt$ at $\ty$ with 
\be
|d\omega_{\tz}(a_{\tz},b_{\tz}) | \geq \mu_0
\ee
where $b_{\tz}$ is the parallel translate of $b$ along the geodesic in $W^u_{\delta_0}(\tz)$ from
$\hz$ to $\tz$, while $a_{\tz}$ is the parallel translate of $a$ along the geodesic in 
$W^s_{\delta_0}(\tz)$ from $\ty$ to $\tz$.}

\bs

\noindent
{\it Remark 3.} In fact, it is clear from the proof of Proposition 6.1 in \cite{kn:St3} that in (ND)
the `parallel translation' in the definition of the vector $b_{\tz}$  can be replaced
by any other uniformly continuous (linear) operator $P_{\hz,\tz} : E^u(\hz) \longrightarrow E^u(\tz)$
($\hz, \tz\in \mt \cap W^u_{\delta_0}(z_0)$). E.g. using a local coordinate system to 
`identify' $E^u(\hz)$ and $E^u(\tz)$ would be good enough. The same applies to the
`parallel translation' in the definition of the vector $a_{\tz}$. In the case of the open billiard 
considered in this paper, (6.1) can be replaced simply by $|d\alpha(a,b)| \geq \mu_0$. 

\bs

As an immediate consequence of the main result\footnote{In fact, the main result in \cite{kn:St3} 
is much more general, however we are not going to discuss it here.} in \cite{kn:St3} 
(see also Sect. 6 there) one gets the following:

\bs

\noindent
{\sc Theorem 6.1.}(\cite{kn:St3}) {\it  Let $\phi_t : M \longrightarrow M$ be a $C^2$ contact flow on a $C^2$ Riemann manifold $M$
and let $\mt$ be a basic set for $\phi_t$ such that the conditions {\rm (P)} and {\rm (ND)} are satisfied for the restriction
of the flow on $\mt$.  Then for any  Lipschitz real-valued function 
$f$ on $\hU$ the Ruelle transfer operators related to  $f$  are eventually contracting}.

\bs

Notice that for open billiard flows both $\wloc^s(x)\cap \Lambda$ and $\wloc^u(x)\cap \Lambda$ are Cantor sets, i.e.
they are infinite compact totally disconnected sets without isolated points. In this particular case we can always 
choose the Markov family $\rr = \{R_i\}_{i=1}^k$ so that the boundary (in $\mt$) of each rectangle $R_i$ is empty and
therefore $U = \hU$.

Next, assume that $K$ is as in Sect. 1. Let $\phi_t$ be the open billiard flow in the exterior of $K$ and let $\mt$
be its non-wandering set. 

The following consequence of Lemma 3.1 shows that under some 
regularity condition, the billiard flow satisfies the condition (ND) on $\mt$.

\bs

\noindent
{\sc Proposition 6.2.}  {\it   Assume that the map $\mt \ni x \mapsto E^u(x)$ is $C^1$. 
Then there exist $z_0\in \mt$, $\delta_0 > 0$  and  $\mu_0 > 0$ such that for any $\delta \in (0,\delta_0]$,
any $z \in \mt \cap W^u_{\delta}(z_0)$ and any unit vector $b \in E^u(z)$ tangent to 
$\mt$ at $z$ there exist $y \in W^s_\delta(z)$ and  a unit vector $a \in E^s(y)$ tangent to $\mt$ at 
$y$ with $|d\alpha(a_z,b) | \geq \mu_0$, where $a_{z}$ is the parallel translate of $a$ along the geodesic in 
$W^s_{\delta_0}(s)$ from $y$ to $z$.}

\bs

\noindent
{\it Proof of Proposition 6.2.} Notice that the standard symplectic form $d\alpha$ in $\R^{2n}$ satisfies
$$|d\alpha(\xi,\eta)| \leq \|\xi\|\, \|\eta\| \quad , \quad  \xi \; , \; \eta \in \R^{2n}\;,$$ 
where we use the standard norm $\|\cdot \|$ in $\R^{2n}$.

 Assume that the map  $\mt\ni x \mapsto E^u(x)$ is $C^1$; then the map  $\mt\ni x \mapsto E^s(x)$ is $C^1$, as well.
Fix $\so$ as in Sect. 3, set $z_0 = \so$, and choosing $ \ep_0 > 0$ and $\delta'_0 \in (0,\ep_0]$ sufficiently small, define
the map 
$$\Phi : W^u_{\delta'_0}(z_0) \longrightarrow W^s_{\ep_0}(z_0)$$ 
as in Sect. 3.
It follows from Lemma 3.1 that there exists a constant $\mu_1 > 0$ (e.g. take $\mu_1 = \kappa/\sqrt{1+C^2}$)
such that $|d\alpha(d\Phi(z_0)\cdot b, b)| \geq \mu_1$ for all unit vectors  $b \in E^u(z_0)$. 
Take $\delta'_0 > 0$ so small that
\be
|d\alpha(d\Phi(z)\cdot b, b)| \geq  \frac{\mu_1}{2} \quad , \quad  z\in W^u_{\delta'_0} (z_0) \; ,\; b \in E^u(z) \;, \;\|b\| = 1\;.
\ee
Further restrictions on $\delta'_0$ will be imposed later.

Next, assuming $\delta'_0 \in (0,\ep_0]$  is sufficiently small, for any $x\in \mt$ and 
$y \in \mt\cap W^u_{\delta'_0}(x)$  the local {\it holonomy map
$\hhu_{x,y} : \mt \cap W^s_{\delta'_0}(x) \longrightarrow \mt \cap W^s_{\ep_0}(y)$
along unstable laminations} is well-defined and uniformly H\"older continuous (see e.g. \cite{kn:HPS} or \cite{kn:PSW}). 
Recall that the map $\hhu_{x,y}$ is defined as follows. Given $x'\in \mt \cap W^s_{\delta'_0}(x)$, there exist a unique
$y'\in W^s_{\ep_0}(y)$  such that $\phi_t(y') \in W^u_{\ep_0}(x')$   for some $t\in \R$, $|t| \leq\ep_0$.
Then we set $\hhu_{x,y}(x') = y'$.  Under the additional condition that the unstable laminations are $C^1$, the 
maps $\hhu_{x,y}$ are $C^1$ as well (see e.g. Fact (2) on  p. 647 in \cite{kn:Ha1}). That is, for each 
$x'\in \mt\cap W^s_{\delta'_0}(x)$ the map $\hhu_{x,y}$ has a {\it linearization} 
$\Lu_{x,y}(x') : E^s(x') \longrightarrow E^s(y')$ 
at $x' \in \mt \cap W^s_{\delta'_0}(x)$ and $\|\Lu_{x,y}(x')\| \leq C_1$ for some constant $C_1 > 0$ 
independent of $x$, $y$ and $x'$.  Notice that $\Lu_{x,y}(x')$ preserves the sets of tangent vectors
to $\mt$, namely if $\xi \in E^s(x')\setminus \{0\}$ is tangent to $\mt$ at $x'$, then 
$\Lu_{x,y}(x') \cdot \xi$  is tangent to $\mt$ at $y'$.


Since the map $\Lu_{z_0,z}(x)$ depends continuously on $z\in W^u_{\delta'_0}(z_0)\cap \mt$ and
$x\in W^s_{\delta'_0}(z_0)\cap \mt$ and $\Lu_{z_0,z_0}(z_0) = I$ (the identity operator), we can 
take $\delta'_0  > 0$  so small that  
\be
\|\Lu_{z_0,z}(x) - I\| \leq \frac{\mu_1}{4}
\ee
for all $z \in \mt\cap W^u_{\delta'_0}(z_0)$ and  $x \in \mt\cap W^s_{\delta'_0}(z_0)$.
We will assume $\delta'_0 > 0$ is chosen so small that for all $z\in W^u_{\delta'_0}(z_0)\cap \mt$,
$y \in W^s_{\delta'_0}(z)\cap \mt$ and unit vectors $a \in E^s(y)$ we have
$\|a - a_z\| \leq \frac{\mu_1}{8 C C_1}$, where $a_z$ is the parallel translate of $a$
along the geodesic on $W^s_{\delta_0}(y)$ from $y$ to $z$.

Finally, take $\delta_0 \in (0,\delta'_0]$ so small that for any $z\in W^u_{\delta_0}(z_0)\cap \mt$
we have $d(\Phi(z), z_0) < \delta'_0$ and $d(z,\hhu_{z_0,z}(\Phi(z))) < \delta'_0$, where $d$ is the
standard distance in $T(\R^n) = \R^{2n}$.
 
Now consider an arbitrary $\delta \in (0,\delta_0]$ and an arbitrary $z\in W^u_\delta(z_0) \cap \mt$.
Let $b \in E^u(z)$ be a unit 
vector tangent to $\mt$ at $z$. Set $x = \Phi(z)$, $a' = d\Phi(z)\cdot b$ and $y = \hhu_{z_0,z}(x)$.
Then $x \in W^s_{\delta'_0}(z_0)\cap \mt$, $a'\in E^s(x)$ is a tangent vector to $\mt$ at
$x$ (see Sect. 3) and $y\in W^s_{\delta'_0}(z)\cap \mt$. Moreover, since $\|b\| = 1$, it follows
from $\|d\Phi\| \leq C$ (see Sect. 3) that $\|a'\| \leq C$. 

Next, the vector $\ta = \Lu_{z_0,z}(x)\cdot a' \in E^s(y)$ is tangent to $\mt$ at $y$
and by (6.3), $\|\ta - a'\| \leq \frac{\mu_1}{4}$.  Moreover, $\|\ta\| \leq C_1 \|a'\| \leq CC_1$.
Hence $a = \frac{\ta}{\|\ta\|} \in E^s(y)$ is a unit vector tangent to $\mt$ at $y$, and 
using (6.2) and (6.3), we get
\begin{eqnarray*}
|d\alpha(a_z,b)| 
& \geq   &  |d\alpha(a,b)| - |d\alpha(a-a_z,b)|
\geq \frac{1}{\|\ta\|} \,|d\alpha(\ta,b)| - \frac{\mu_1}{8 C C_1}\\
& \geq   & \frac{1}{C C_1}\, \left[ |d\alpha(a',b)| - |d\alpha(\ta-a',b)|\right] - \frac{\mu_1}{8 C C_1}\\
& \geq   & \frac{1}{C C_1}\, \left[ |d\alpha(d\Phi(z)\cdot b,b)| - \frac{\mu_1}{4}\right] 
- \frac{\mu_1}{8 C C_1} \geq \frac{\mu_1}{4CC_1} - \frac{\mu_1}{8 C C_1} = \mu_0 \;,
\end{eqnarray*}
where $\mu_0 = \frac{\mu_1}{8 C C_1}$.
This proves the assertion.
\endofproof

\bs

From Theorem 6.1 and Proposition 6.2 one derives the following.

\bs

\noindent
{\sc Theorem 6.3.} {\it  Assume that the billiard flow $\phi_t : \mt \longrightarrow \mt$ 
satisfies the condition {\rm (P)}  on its non-wandering set $\mt$. Then for any  Lipschitz real-valued function 
$f$ on $U$ the Ruelle transfer operators related to  $f$  are eventually contracting}.

\bs

\noindent
{\it Proof of Theorem 6.3.}  As mentioned in Sect. 1, the condition (P) implies 
that  the map\\  $\mt\ni x \mapsto E^u(x)$ is $C^1$. Then by Proposition 6.2, $\phi_t$ satisfies the
condition {\rm (ND)} on $\mt$. Now applying Theorem 6.1 proves the assertion.
\endofproof

\bs

Results of this kind were first established  by Dolgopyat (\cite{kn:D}) for some Anosov flows 
(i.e. $\mt = M$, a compact Riemann manifold). His results apply to geodesic flows on any compact surface
(for any $f$), and to transitive Anosov flows on compact Riemann manifolds with $C^1$ jointly 
non-integrable local stable and unstable foliations  for the Sinai-Bowen-Ruelle potential  $f = \log \det (d\phi_\tau)_{|E^u}$.

As one can see, Theorems 6.1  and 6.3 work for any potential. Theorem 6.3
generalizes the result in \cite{kn:St2} which deals with open billiard flows in the plane.


It should be mentioned that Dolgopyat type estimates for pinched open billiard flows have already been used 
in \cite{kn:PS1}, \cite{kn:PS2} and \cite{kn:PS3} to obtain some rather non-trivial results. 
The main result  in \cite{kn:PS1} provides 
existence of an analytic continuations of the cut-off resolvent of the 
Dirichlet Laplacian  in $\R^n\setminus K$ in a 
horizontal strip above the level of absolute convergence and polynomial estimates for the norm of the cut-off
resolvent in  such a domain. The Dolgopyat type estimates for the open billiard 
flow in $\R^n\setminus K$ play a
significant role in the proof. These estimates are also essential for the proof 
of the main result in \cite{kn:PS2}
which deals with estimates of correlations for pairs of closed billiard 
trajectories for open billiards. Previous results
of this kind were established in \cite{kn:PoS2} for geodesic flows on surfaces of negative curvature.
Finally, in a very recent preprint \cite{kn:PS3}, using  Theorem 6.3 
a fine asymptotic was obtained  for the number of closed billiard trajectories in $\mt$
with primitive periods lying in exponentially shrinking intervals 
$(x - e^{-\delta x}, x + e^{-\delta x})$, $\delta > 0$, $x \to + \infty.$

As in \cite{kn:St3}, using  Theorem 6.3  and an argument of Pollicott and Sharp \cite{kn:PoS1}, 
we get some rather significant consequences about the {\it Ruelle zeta function}
$\zeta(s) = \prod_{\gamma} (1- e^{-s\ell(\gamma)})^{-1}\;.$
Here $\gamma$ runs over the set of primitive  closed orbits of $\phi_t: \mt \longrightarrow \mt$
and $\ell(\gamma)$ is the least period of $\gamma$.  Let $h_T$ be the
 {\it topological entropy} of $\phi_t$ on $\mt$.\\

\noindent
{\sc Corollary 6.4.} {\it Under the assumptions in Theorem 6.3, the  zeta function $\zeta(s)$
of the flow $\phi_t: \mt \longrightarrow \mt$ has an analytic  and non-vanishing continuation in a 
half-plane $\Re(s) > c_0$ for some $c_0 < h_T$ except for a simple pole at $s = h_T$.  Moreover, there exists
$c \in (0, h_T)$ such that
$$\pi(\lambda) = \# \{ \gamma : \ell(\gamma) \leq \lambda\} = \li(e^{h_T \lambda}) + O(e^{c\lambda})$$
as $\lambda\to \infty$, where $\di \li(x) = \int_2^x \frac{du}{\log u} \sim \frac{x}{\log x}$ as 
$x \to \infty$. }

\bigskip

As another consequence of Theorem 6.3  and the procedure described in \cite{kn:D} 
one gets exponential decay of
correlations for the  flow $\phi_t : \mt \longrightarrow \mt$. 

Given $\alpha > 0$ denote by $\ff_\alpha(\mt)$ the set of {\it H\"older 
continuous functions} with H\"older exponent $\alpha$ and by $\| h\|_\alpha$ the 
{\it H\"older constant} of $h\in \ff_\alpha(\mt)$. \\

\noindent
{\sc Corollary 6.5.} {\it   Under the assumptions in Theorem 6.3,  let $F$ be a H\"older continuous 
function on $\mt$  and let 
$\nu_F$ be the Gibbs measure determined by $F$ on $\mt$. Assume in addition that the boundary of 
$K$ is at least $C^5$.
Then for every $\alpha > 0$ there exist constants $C = C(\alpha) > 0$ and $c = c(\alpha) > 0$ such that 
$$\left| \int_{\mt} A(x) B(\phi_t(x))\; d\nu_F(x) - 
\left( \int_{\mt} A(x)\; d\nu_F(x)\right)\left(\int_{\mt} B(x) \; d\nu_F(x)\right)\right|
\leq C e^{-ct} \|A\|_\alpha \; \|B\|_\alpha \;$$
for any two functions $A, B\in \ff_\alpha(\mt)$.}\\

One would expect that much stronger results could be established by using the techniques recently developed in  
\cite{kn:BKL}, \cite{kn:L}, \cite{kn:BG},  \cite{kn:GL}, \cite{kn:T} (see the references there, as well). Still, there are not very many
results of this kind. In fact, for dimensions higher than two the author is not aware  of any other results of this kind concerning billiard flows.
What concerns billiards in general, bounds of correlation decay known so far concern mostly the corresponding discrete dynamical system 
(generated by the  billiard ball map from boundary to boundary) -- see  \cite{kn:BSC}, \cite{kn:Y} and \cite{kn:Ch2}.  See also \cite{kn:ChZ}
and the references there for some related results. Recently, a sub-exponential decay of correlations for Sinai billiards in the plane was established 
by Chernov (\cite{kn:Ch3}).
For open billiard flows in the plane exponential decay of correlations was proved in \cite{kn:St2} (as a consequence of the
Dolgopyat type estimates established there).


\end{document}